\documentclass []{article}
\usepackage{amssymb}
\def\qmod#1#2{{\hbox{}^{\displaystyle{#1}}}\!\big/\!\hbox{}_{
\displaystyle{#2}}}

\def\resto#1#2{{
#1\hskip 0.4ex\vline_{\hskip 0.4ex\raisebox{-1ex}
{{${\scriptstyle #2}$}}}}}

\newfam\msbfam

\font\tenmsb=msbm10

\font\sevenmsb=msbm10 at 7pt
\font\fivemsb=msbm10 at 5pt

\textfont\msbfam=\tenmsb
\scriptfont\msbfam=\sevenmsb
\scriptscriptfont\msbfam=\fivemsb

\def\Bbb{\fam\msbfam\tenmsb}

\def\C{{\Bbb C}}

\def\H{{\Bbb H}}

\def\N{{\Bbb N}}

\def\Q{{\Bbb Q}}
\def\R{{\Bbb R}}

\def\Z{{\Bbb Z}}

\def\qed {\hfill\vrule height6pt width6pt depth0pt \bigskip}

\def\map{\longrightarrow}
\def\textmap#1{\mathop{\vbox{\ialign{
                                 ##\crcr
     ${\scriptstyle\hfil\;\;#1\;\;\hfil}$\crcr
     \noalign{\kern 0pt\nointerlineskip}
     \rightarrowfill\crcr}}\;}}

\def\textlmap#1{\mathop{\vbox{\ialign{
                                 ##\crcr
     ${\scriptstyle\hfil\;\;#1\;\;\hfil}$\crcr
     \noalign{\kern-1pt\nointerlineskip}
     \leftarrowfill\crcr}}\;}}

\newfam\meuffam

\font\tenmeuf=eufm10
\font\sevenmeuf=eufm7
\font\fivemeuf=eufm5

\textfont\meuffam=\tenmeuf
\scriptfont\meuffam=\sevenmeuf
\scriptscriptfont\meuffam=\fivemeuf

\def\germ{\fam\meuffam\tenmeuf}

\def\g{{\germ g}}

\newtheorem{sz}{Satz}[section]
\newtheorem{thry}[sz]{Theorem}
\newtheorem{pr}[sz]{Proposition}
\newtheorem{re}[sz]{Remark}
\newtheorem{co}[sz]{Corollary}
\newtheorem{dt}[sz]{Definition}
\newtheorem{lm}[sz]{Lemma}

\begin{document}
\def\Pr{{\rm Pr}}
\def\tr{{\rm Tr}}
\def\End{{\rm End}}
\def\Aut{{\rm Aut}}
\def\Spin{{\rm Spin}}
\def\U{{\rm U}}
\def\SU{{\rm SU}}
\def\SO{{\rm SO}}
\def\PU{{\rm PU}}
\def\GL{{\rm GL}}
\def\spin{{\rm spin}}
\def\su{{\rm su}}
\def\so{{\rm so}}
\def\ub{\underbar}
\def\pu{{\rm pu}}
\def\Pic{{\rm Pic}}
\def\Iso{{\rm Iso}}
\def\NS{{\rm NS}}
\def\deg{{\rm deg}}
\def\Hom{{\rm Hom}}
\def\Aut{{\rm Aut}}
\def\h{{\germ h}}
\def\Herm{{\rm Herm}}
\def\Vol{{\rm Vol}}
\def\pf{{\bf Proof: }}
\def\id{{\rm id}}
\def\i{{\germ i}}
\def\im{{\rm im}}
\def\rk{{\rm rk}}
\def\ad{{\rm ad}}
\def\h{{\bf H}}
\def\coker{{\rm coker}}
\def\dbar{\bar{\partial}}
\def\Lo{{\Lambda_g}}
\def\niq{=\kern-.18cm /\kern.08cm}
\def\Ad{{\rm Ad}}
\def\RSU{\R SU}
\def\ad{{\rm ad}}
\def\dva{\bar\partial_A}
\def\da{\partial_A}
\def\p{\partial\bar\partial}
\def\sp{\Sigma^{+}}
\def\sm{\Sigma^{-}}
\def\spm{\Sigma^{\pm}}
\def\smp{\Sigma^{\mp}}
\def\Tors{{\rm Tors}}
\def\st{{\rm st}}
\def\s{{\rm s}}
\def\oo{{\scriptstyle{\cal O}}}
\def\ooo{{\scriptscriptstyle{\cal O}}}
\def\sw{Seiberg-Witten }
\def\pa{\partial_A\bar\partial_A}
\def\Dr{{\raisebox{0.15ex}{$\not$}}{\hskip -1pt {D}}}
\def\gr{{\scriptscriptstyle|}\hskip -4pt{\g}}
\def\subsetint{{\  {\subset}\hskip -2.45mm{\raisebox{.28ex}
{$\scriptscriptstyle\subset$}}\ }}
\def\ra{\rightarrow}
\def\kod{{\rm kod}}

\title{The pseudo-effective cone of  a non-K\"ahlerian surface and applications}
\author{Andrei Teleman\\
teleman@cmi.univ-mrs.fr}
\maketitle


\begin{abstract}
We describe the positive cone and the pseudo-effective cone of a non-K\"ahlerian surface. We use 
these results for two types of applications: 
\begin{enumerate}
\item Describe the set $\sigma(X)\subset \R$ of possible total Ricci scalars associated with
Gauduchon metrics of fixed volume 1 
on a fixed non-K\"ahlerian surface, and  decide whether
the assignment
$X\mapsto\sigma(X)$ is a deformation invariant.
\item Study the stability   of the canonical extension $$0\to {\cal K}_X\to {\cal
A}\to{\cal O}_X\to 0$$ of a class VII surface $X$ with positive $b_2$.
This extension plays an important role in our strategy to prove the GSS conjecture using gauge theoretical methods
\cite{Te2}, \cite{Te3}.
\end{enumerate}

Our main tools are Buchdahl's ampleness criterion for non-K\"ahlerian surfaces \cite{Bu2} and
the recent results of Dloussky-Oeljeklaus-Toma \cite{DOT} and Dloussky \cite{D} 
  on class VII surfaces with curves.
\end{abstract}



\section{Introduction}

In this paper we study certain complex geometric and differential geometric properties   of
non-K\"ahlerian surfaces. The first problem  we treat   is the following:\\

{\it  Describe explicitly the   pseudo-effective  cone of a non-K\"ahlerian surface and
compare it with the effective cone.}\\

By definition, the pseudo-effective  cone of a complex surface is the subset of the Bott-Chern
cohomology space $H^{1,1}_{BC}(X,\R)$ consisting of $i\bar\partial\partial$-classes which are
represented by closed positive currents. The effective cone is just the cone generated by classes
associated with effective divisors.
 
We will solve completely this problem showing that the pseudo-effective cone is determined in a
simple way by the finite set of irreducible effective divisors with negative self-intersection.
The proof is  based on a version of Buchdahl's ampleness criterion \cite{Bu2}, which will be
explained in the first section. This criterion will provide a simple description of the {\it
positive cone} of a non-K\"ahlerian surface, i.e. the cone of $i\bar\partial\partial$-closed
(1,1)-classes (modulo
$i\bar\partial\partial$-exact forms)  associated with Gauduchon metrics \cite{G}.

We point out that {\it all our results do not make use of the GSS conjecture; in particular
they hold for  the still non-classified class $VII$-surfaces with   second Betti number
$b_2>1$.}\\

Our description of the pseudo-effective cone will allow us to solve the following two problems.

\begin{enumerate}
\item {\it Determine the possible values of the total {\it Ricci scalars} of the
Gauduchon metrics with volume 1 on a given non-K\"ahlerian surface. }

For a Hermitian metric $g$ we denote by $s_g$ the {\it  Ricci scalar} of $g$
(see \cite{G}), which is defined by the formula
$$s_g:=i\Lambda_g \tr( F_{A_g})\ ,
$$
where $A_g$ is the Chern connection associated with $g$ and the holomorphic structure on the
tangent bundle. In the non-K\"ahlerian framework, $s_g$ does not coincide in general with the
scalar curvature of the Riemannian metric $g$.  The total Ricci scalar of $g$ is  
defined by  
$$\sigma_g:=\int_X s_g vol_g=\int_X i\tr(F_{A_g})\wedge \omega_g\ .
$$
Let ${\cal G}(X)$ be the set of Gauduchon metrics of $X$. Our problem is to determine the set
$$\sigma(X):=\{\sigma_g|\ g\in{\cal G}(X),\ \int_X vol_g=1\}.
$$
We will see for instance that, for a (blown up) Inoue surface $X$, one has
$\sigma(X)=(-\infty,0)$, which might be surprising for a surface with $\kod(X)=-\infty$. We
will also answer the following natural question:\\

{\it Is the assignment   $X\mapsto \sigma(X)\subset\R$ a deformation invariant ?}\\

Using certain families of class $VII$-surfaces, we will see that the answer is in general
negative.

\item {\it The stability of the canonical extension of a class $VII$-surface.}

Let $X$ be a class $VII$-surface (i.e. a surface with $b_1(X)=1$ and $\kod(X)=-\infty$). Such a
surface has $h^1(X,{\cal O}_X)=1$  so, by Serre duality, one also  has $h^1(X,{\cal
K}_X)=1$.  The {\it canonical extension} of $X$ is defined to be the unique (up to the
natural
$\C^*$-action on
${\rm Ext}^1({\cal O}_X,{\cal K}_X)$) nontrivial extension of the form
\begin{equation} \label{canext}
0\map {\cal K}_X\map {\cal A}\map {\cal O}_X\map 0
\end{equation}

The question here is whether there exists Gauduchon metrics on $X$ with respect to which ${\cal
A}$ is stable. 

  We will see that, excepting certain very special surfaces with global spherical
shell, every minimal class VII surface with positive $b_2$ admits Gauduchon metrics
$g$ for which the bundle ${\cal A}$  is
$g$-stable. The motivation for this problem is the following:

For any   topologically trivial line bundle ${\cal
L}\in\Pic^0(X)\simeq\C^*$ with ${\cal L}^{\otimes 2}\not\simeq{\cal O}_X$ one has ${\rm
Ext}^1({\cal L},{\cal K}_X\otimes{\cal L}^{-1})=0$, so there are no non-trivial
extensions of ${\cal L}$ by ${\cal K}_X\otimes{\cal L}^{-1}$.  On the other hand,
the dimension of the moduli space
${\cal M}^{\rm st}(0,{\cal K}_X)$ of stable rank 2-bundles ${\cal E}$ with $c_2({\cal E})=0$
and
$\det({\cal E})={\cal K}_X$ is
$b_2(X)$.  

Therefore, although the extension (\ref{canext}) is rigid (it  cannot be deformed
in another extension of the form $0\to{\cal K}_X\otimes{\cal
L}^{-1}\to{\cal E}\to{\cal L}\to 0$ with $[{\cal L}]\in\Pic^0(X)$), its  central
term
${\cal A}$ cannot be rigid for $b_2(X)>0$.

As in \cite{Te2}, using the Kobayashi-Hitchin correspondence on non-K\"ahlerian surfaces
(\cite{Bu1}, \cite{LY}, \cite{LT}) one can prove that, if
$X$ had no curve and
$b_2(X)\leq 3$, the connected component of $[{\cal A}]$ in ${\cal M}^{\rm st}(0,{\cal K}_X)$
would be a   smooth  compact manifold containing both filtrable and non-filtrable 
points. This is the starting point of our strategy to prove the GSS conjecture
using gauge theoretical methods.

\end{enumerate}

\section{Buchdahl's ampleness criterion and positiveness}

In \cite{Bu2} Buchdahl proved an interesting ampleness criterion for
(non-algebraic) complex surfaces; surprisingly, his statement is very much
similar to the algebraic geometric Nakai-Moishezon criterion. This result 
suggests that certain fundamental purely algebraic geometric theorems might
have natural extensions to the non-algebraic and even  non-K\"ahlerian
framework; the difficulty is   to find the correct complex geometric analogues of
the algebraic geometric notions involved in the original statement.

\begin{thry}\label{ampl1} \cite{Bu2} Let $X$ be a compact complex surface equipped
with a positive $\bar\partial\partial$-closed (1,1)-form $\omega_0$ and
let $\varphi$ be a smooth real $\bar\partial\partial$-closed (1,1)-form
satisfying
\begin{enumerate}
\item $\int_X \varphi\wedge \varphi>0,$
\item $\int_X \varphi\wedge\omega_0>0,$
\item  $\int_D \varphi >0$ for every irreducible effective divisor with $D^2<0$.
\end{enumerate}
Then there is a smooth function $\psi$ on $X$ such that
$\varphi+i\bar\partial\partial\psi$ is positive.
\end{thry}

This result shows that it is very natural to extend the fundamental
algebraic geometric notion ``positive cone`" to the non-algebraic
non-K\"ahlerian framework in the following way:

Set
$$Q_0:=i\bar\partial\partial:A^0(X,\R)\map A^{1,1}(X,\R)\ ,\
Q_1:=i\bar\partial\partial:A^{1,1}(X,\R)\map A^{2,2}(X,\R) \ ,$$
$${\cal
H}(X):=\qmod{\ker(Q_1)}{\im(Q_0)}\ .
$$

It is easy to see that ${\im(Q_0)}$ is closed: it suffices to choose a Hermitian metric $g$ on $X$ and to 
note that the operator $P_g:=\Lambda_g\circ Q_0$ is elliptic. Therefore
${\cal H}(X)$ is a Fr\'echet space. It is not difficult to see that this space is
infinite dimensional; it contains the finite dimensional Bott-Chern
cohomology space (see \cite{BHPV}, p. 148)
$$H_{BC}^{1,1}(X,\R):=\qmod{\ker(d:A^{1,1} (X,\R)\map
A^3 (X,\R))}{i\bar\partial\partial(A^0(X,\R))}\ .
$$
\begin{dt} Let $X$ be a compact complex surface and ${\cal G}(X)$ the space of Gauduchon metrics
on $X$. The positive cone of
$X$ is the open cone   ${\cal H}_+(X)\subset {\cal H}(X)$ defined by   
$${\cal H}_+(X):=\{[\omega_g]|\ g\in{\cal G}(X)\}\ .$$ 
\end{dt}

Note  that one has a natural well defined intersection form
$${\cal H}(X)\times {\cal H}(X)\to \R
$$
given by $[\eta]\cdot[\mu]\mapsto \int_X \eta\wedge\mu$. Moreover, every
real $i\bar\partial\partial$-closed   $(1,1)$-current  $u$ defines a
linear form $\langle\cdot, u\rangle:{\cal H}(X)\to\R$.

Buchdahl's criterion says that the elements $h$ of the positive cone ${\cal H}_+$
are characterized by  the system of  inequalities:
\begin{enumerate}
\item $h^2>0$,
\item $h\cdot[\omega_0]>0$,
\item $\langle h, [D]\rangle >0$ for every irreducible effective 
divisor $D$ with $D^2<0$.
\end{enumerate}

In the non-K\"ahlerian case ($b_1(X)$ odd), one can reformulate  this
criterion, by replacing the class $[\omega_0]$ in the second inequality  
with the class of an {\it exact} form. This modification, which is explained in
detail  below, is very useful, because  all linear inequalities in the resulting
system will be associated with classes in the Bott-Chern cohomology
space $H_{BC}^{1,1}(X,\R)$.\\

For a complex surface $X$ we put

$$B^{1,1}(X,\R):=d(A^1(X,\R))\cap A^{1,1} (X,\R)\supset
i\bar\partial\partial(A^0 (X,\R))\ ,
$$
$$ H^{1,1}(X,\R):=\qmod{\ker(d:A^{1,1} (X,\R)\map
A^3 (X,\R))}{B^{1,1}(X,\R)}\subset H^2(X,\R)\ .
$$

Some of the statements in the following lemma are probably known. We include
short proofs for completeness.

\begin{lm} Let $g$ be a Gauduchon metric  on $X$.  
\begin{enumerate} 
\item One has the exact sequences
$$0\map \qmod{B^{1,1} (X,\R)}{i\bar\partial\partial(A^0(X,\R))}\map
H_{BC}^{1,1}(X,\R)\map H^{1,1} (X,\R)\map 0
$$
$$
0\map i\bar\partial\partial(A^0(X,\R)) \map B^{1,1} (X,\R)
\textmap{\langle
\cdot\ ,\ \omega_g\rangle }\R
$$
\item The canonical map $H^1(X,i\R)\to H^1(X,{\cal O})$ is injective.
\item The following conditions are equivalent
\begin{enumerate} 
\item $\langle\cdot,\omega_g\rangle$ vanishes identically on
$B^{1,1}(X,\R)$.
\item $H_{BC}^{1,1}(X,\R)= H^{1,1} (X,\R)$
\item The  natural monomorphism $H^1(X,i\R)\to H^1(X,{\cal O}_X)$ is
surjective.
\item $b_1(X)$ is even.
\end{enumerate}
\item When $b_1(X)$ is odd, one has an exact sequence
$$0\map \Gamma(X)\map
H_{BC}^{1,1}(X,\R)\map H^{1,1}(X,\R)\map 0
$$
 where $\Gamma:=\qmod{B^{1,1}(X,\R)}{i\bar\partial\partial(A^0(X,\R))}$
is a line which is identified with $\R$ via
$\langle\cdot,\omega_g\rangle$.
\end{enumerate}
\end{lm}
\pf 
1. The  first exact sequence is obvious. For the second, let $\alpha\in
B^{1,1}(X,\R)$ such that
$$\int_X \alpha\wedge \omega=0
$$
The Laplace equation $i\Lambda_g \bar\partial\partial u=\Lambda_g\alpha=0$
is solvable, because the image of the elliptic operator
$P_g:=i\Lambda_g \bar\partial\partial$ is precisely $\ker \langle \cdot
,\omega_g\rangle$; let $u_0$ be a solution of this equation. The form 
$i\bar\partial\partial u_0-\alpha$ is exact and anti-selfdual, so it
vanishes. \\
\\
2. Let $\alpha=-a^{1,0}+a^{0,1}$ be a closed imaginary form (where
$a^{0,1}=\overline{a^{1,0}}$) such that $a^{0,1}$ is
$\bar\partial$-exact, and let $u\in A^0(X,\C)$ such that $\bar\partial
u=a^{0,1}$. This implies $\partial \bar u= a^{1,0}$. Since $\alpha$ is
closed, we get
$$ \partial\bar\partial u-  \bar\partial\partial \bar
u=-2 \bar\partial\partial ( {\rm Re} (u) ) =0\ ,
$$
hence, we can suppose that $u$ is purely imaginary. Then  one gets
immediately $(\partial+\bar\partial) u=\alpha$. 
\\

3. By the second exact sequence in 1.,  the statements $(a)$ and $(b)$ are
both equivalent to the equality
\begin{equation}\label{id}
\bar\partial\partial(A^0(X,\R))= B^{1,1} (X,\R)
\end{equation}

To prove $(b)\Rightarrow (c)$ it suffices to show that (\ref{id}) implies
the surjectivity of the map $H^1(X,i\R)\to H^1(X,{\cal O}_X)$.
Let $[\beta^{0,1}]\in H^1(X,{\cal O}_X)$. It suffices to find $\varphi\in
A^0(X,\C)$ such that the form 
$\alpha^{0,1}:=\beta^{0,1}+\bar\partial\varphi$ satisfies
$$ \partial\alpha^{0,1}-\bar\partial \alpha^{1,0}=0
$$
(where $\alpha^{1,0}:=\bar \alpha^{0,1}$). Indeed, if we find such a function $\varphi$, the de
Rham class of $\alpha:=\alpha^{0,1}-\alpha^{1,0}$ will be a preimage of the Dolbeault class of
$\beta$ under the natural morphism $H^1(X,i\R)\to H^1(X,{\cal O}_X)$.
 This equation becomes

$$\partial(\beta^{0,1}+\bar\partial\varphi)-\bar\partial (\beta^{1,0}+  \partial\bar\varphi)=0\ .
$$
We will show that there exists a real solution of this equation. For a real function $\varphi$ the
equation becomes 
$$2i\partial\bar\partial\varphi=i(\bar\partial  \beta^{1,0}- \partial \beta^{0,1})
$$

The form  $i(\bar\partial  \beta^{1,0}- \partial \beta^{0,1})$ can be written as $d(i
\beta^{1,0}- i
\beta^{0,1})$  (because $\bar\partial \beta^{1,0}=\partial\beta^{10}=0$), hence it is an
exact $(1,1)$ form. Therefore, by hypothesis, it belongs to the image of the operator
$i\partial\bar\partial$.  The implication $(c)\Rightarrow (d)$ is obvious.  

The implication $(d)\Rightarrow (b)$ is well known: it follows from  Theorems 2.8, 2.10 and
Corollary 13.8 in \cite{BHPV}.
\\ \\
4. This follows directly from 1 and 3.
\qed
\begin{re} Let $X$ be a complex surface with $b_1(X)$ odd.
Since the space of Gauduchon metrics on $X$ is connected,  
the orientation of the ``exact line''
$\Gamma(X)\subset  
H_{BC}^{1,1}(X,\R)$ induced by $\langle\cdot,\omega_g\rangle$ is well defined.  Let $\gamma_0$  be
a positive generator of this line.  One has
\begin{equation}
\langle \gamma_0 \cdot [\omega_g]\rangle >0
\end{equation}
for {\it every} Gauduchon  metric $g$ on $X$.
\end{re}

Recall now that the Bott-Chern cohomology space $H_{BC}^{1,1}(X,\R)$ can be also introduced using
currents (see
\cite{BHPV}, p. 148):
$$H_{BC}^{1,1}(X,\R):=\qmod{\ker(d:A^{1,1} (X,\R)\map
A^3 (X,\R))}{i\bar\partial\partial(A^0(X,\R))}=$$
$$\qmod{\ker(d:{\cal D}'_{1,1} (X,\R)\map
{\cal D}'_{0} (X,\R))}{i\bar\partial\partial({\cal D}'_{2,2} (X,\R))}\ .
$$
\begin{pr}\label{exactpos} Let $X$ be a complex surface with odd $b_1(X)$. Then the
Bott-Chern cohomology class $\gamma_0$ is represented by an exact
positive current.
\end{pr}
\pf It is well known that a complex surface with odd first Betti number admits a non-trivial
exact positive (1,1)-current (see \cite{Bu2}, \cite{La}). This is a refinement -- valid for
surfaces -- of the general  Harvey-Lawson's characterization of non-K\"ahlerianity \cite{HL}.

Let $v$ be non-trivial
exact positive (1,1)-current on $X$. Then one has obviously $\langle v,\omega_g\rangle>0$, for
every Gauduchon metric on $X$, so the Bott-Chern cohomology class $[v]$ is a positive multiple of
the generator $\gamma_0$.
\qed

We claim that, when $b_1(X)$ is odd, Buchdahl's ampleness criterion is equivalent
to the following:
\begin{thry}\label{ampl2} Let $X$ be a surface with $b_1(X)$ odd. The elements $h$
of the positive cone
${\cal H}_+(X)$ are characterized by the inequalities
\begin{enumerate}
\item $h^2>0$,
\item $h\cdot \gamma_0>0$,
\item $\langle h,D\rangle>0$ for every irreducible effective divisor with $D^2<0$.
\end{enumerate}
\end{thry} 
\pf Let $g_0$ be a  Gauduchon metric $g_0$ on $X$ and $\omega_0$ the corresponding
form. For any $t\geq 0$, the class $[\omega_0]+t\gamma_0$ satisfies the three inequalities
in Buchdahl's criterion, because
$$([\omega_0]+t\gamma_0)^2=[\omega_0]^2+2t\gamma_0\cdot [\omega_0] \ ,\ 
([\omega_0]+t\gamma_0)\cdot[\omega_0]=[\omega_0]^2+ t\gamma_0\cdot [\omega_0]\ ,$$
$$
\langle [\omega_0]+t\gamma_0 ,D\rangle=\langle [\omega_0],D\rangle\ ,
$$
for every effective divisor $D$.  Therefore $[\omega_0]+t\gamma_0$ is still the class of the
K\"ahler form, say $\omega_t$, of a Gauduchon metric $g_t$ on $X$.

Let $h\in {\cal H}$ be a class satisfying the three inequalities in the hypothesis.
For sufficiently large $t>0$, one has $h\cdot ([\omega_0]+t\gamma_0)>0$. Therefore $h$
satisfies   Buchdahl's original criterion (for $\omega_t$ instead of $\omega_0$).
\qed
\section{Effectiveness and pseudo-effectiveness}

Let ${\cal D}'_{1,1}(X,\R)$, ${\cal D}'_{0}(X,\R)$, ${\cal
D}'_{2,2}(X,\R)$ be the dual spaces of
$A^{1,1}(X,\R)$,  $A^{0}(X,\R)$, $A^{2,2}(X,\R)$ respectively.
Using again the ellipticity of the operator $P_g=i\Lambda_g\bar\partial\partial$ associated with a
Hermitian metric, one gets easily

\begin{re} The image $Q_1^*({\cal D}'_{2,2}(X,\R))\subset {\cal
D}'_{1,1}(X,\R)$ is closed in
${\cal D}'_{1,1}(X,\R)$ with respect to the weak topology.
\end{re}

\begin{co}\label{dual} The dual space of the quotient $\qmod{{\cal D}'_{1,1}(X)}{Q_1^*({\cal
D}'_{2,2}(X))}$ is naturally isomorphic to $
\ker Q_1$.
\end{co}
\pf Indeed, by a well-known result about duality in locally convex spaces, the dual
space of $\qmod{{\cal D}'_{1,1}(X)}{Q_1^*({\cal
D}'_{2,2}(X))}$ is just the subspace of ${\cal D}'_{1,1}(X)^*=A^{1,1}(X)$ consisting
of functionals which vanish on $Q_1^*({\cal
D}'_{2,2}(X))$. It's easy to see that this subspace is just 
$\ker Q_1$.
\qed

\begin{re} Let $X$ be complex surface with $b_1(X)$ odd. The set ${\cal
D}ou(X)_-^{\rm irr}$ of irreducible effective divisors with negative
self-intersection is finite.
\end{re}
\pf
Indeed, a complex surface $X$ with odd $b_1(X)$ is either an elliptic fibration, or a class VII
surface. In the first case, $X$ cannot contain horizontal divisors (because otherwise $X$ would be
algebraic) and, on the other hand, the generic fibre has vanishing self-intersection. Therefore, the
irreducible effective divisors with negative self-intersection must be components of the (finitely
many) singular fibres.

If $X$ is not an elliptic fibration, it must be a class VII surface of vanishing
algebraic dimension, so it has only finitely many  irreducible
effective divisors.
\qed

For a finite subset $A\subset E$ of a real vector space $E$, we denote by $[A]$ the convex hull
of $A$ and by $CA$ the cone over $A$ 
$$CA:=\left\{\sum_{a\in A} t_a a|\ t_a\geq 0\ \forall a\in A\right\}\ .
$$

\begin{lm}\label{separation} Let $E$ be a locally convex space, and $A$, $B\subset  E$  non-empty
finite sets, such that
$0\not\in [A]\cup [B]$. Then the following conditions are equivalent:
\begin{enumerate}
\item $CA\cap CB=\{0\}$,
\item There exists a continuous linear form $u\in E^*$ such that 
$$\resto{u}{A}<0\hbox{ and }\resto{u}{B}>0\ .$$
\end{enumerate}
\end{lm}
\pf The implication $2.\Rightarrow 1.$ is obvious. For the other implication, apply the standard
separation theorem in locally convex spaces to the closed convex set $CA$ and the compact
convex set
$[B]$. These sets are disjoint  because, since $0\not\in [B]$,  any intersection
point would be a non-zero element in
$CA\cap CB$. We get a continuous  linear form 
$v\in E^*$ such that
$$\sup_{x\in CA} v(x) < \inf_{y\in [B]} v(y)
$$
Since $0\in CA$, we get $\sup\limits_{x\in CA} v(x)\geq 0$. On the other hand one must have
$v(a)\leq 0$ for every $a\in A$  because, otherwise, one would obviously have  
$$\sup_{x\in CA} v(x)=+\infty\ .$$
Therefore $\sup\limits_{x\in CA} v(x)=0$ and $\inf\limits_{y\in [B]} v(y)>0$. This means
$$
\resto{v}{A}\leq 0\hbox{ and }\resto{v}{B}>0\ .
$$

Since $0\not\in [A]$ and $[A]$ is compact convex set, there exists a continuous linear functional
$w\in E^*$ such that
$\resto{w}{A}>0$. Setting  $u:=v+\varepsilon w$ for sufficiently small
$\varepsilon>0$, we get the desired inequalities. 
\qed
\begin{thry}\label{psef} 
Let $X$ be complex surface with $b_1(X)$ odd.  Let
$c\in H^{1,1}_{BC}(X,\R)$ be a Chern-Bott cohomology class such that
$\langle c,\omega_g\rangle\geq0$ for every Gauduchon metric $g$ on $X$. Then
  $c\in C(\{{\gamma_0}\}\cup {\cal
D}ou(X)_-^{\rm irr})$.
\end{thry}
\pf   Suppose $c\ne 0$. If the claim was not true, then by Lemma \ref{separation}   one would find
a closed linear hyperplane separating
$c$ from this cone.  Therefore, by Corollary \ref{dual}, there would exist a smooth
$i\bar\partial\partial$-closed
$(1,1)$-form
$\eta\in A^{1,1}(X)$ such that 
\begin{enumerate}
\item $\langle \eta,\gamma_0\rangle>0$ and $\langle \eta,D\rangle>0$ for all 
$D\in{\cal
D}ou(X)_-^{\rm irr}$
\item  $\langle \eta, c\rangle <0$.
\end{enumerate}

One has $([\eta]+ t \gamma_0)^2=[\eta]^2+2t[\eta]\cdot \gamma_0$, which becomes positive
for sufficiently large $t$. Therefore $ [\eta]+ t \gamma_0$ satisfies the three
assumptions of the ampleness criterion  given by  Theorem  \ref{ampl2}. This would
give a Gauduchon metric $g$ on $X$  with
$$\langle \omega_g,c\rangle=\langle [\eta]+ t \gamma_0, c\rangle=\langle [\eta] ,
c\rangle<0\ ,
$$
which contradicts the hypothesis.
\qed

Putting together Theorem \ref{psef} and Proposition \ref{exactpos}, we get 
\begin{co}\label{equiv} Let $X$ be a surface satisfying the assumptions of Theorem  \ref{psef},
and let $u$ be  a real, closed
$(1,1)$-current. Then the following conditions are equivalent.
\begin{enumerate}
\item  $\langle u,\omega_g\rangle\geq0$ for every Gauduchon metric $g$ on $X$.
\item The Bott-Chern cohomology class $[u]\in H^{1,1}_{BC}(X,\R)$ of  $u$  belongs to the cone 
$C(\{{\gamma_0}\}\cup {\cal D}ou(X)_-^{\rm irr})$

\item The Bott-Chern cohomology class $[u]\in H^{1,1}_{BC}(X,\R)$ of  $u$ is represented
by a closed positive current.

\end{enumerate}
\end{co}

Following the standard terminology used in the algebraic and K\"ahlerian
case \cite{De}, we define

\begin{dt} 
\item A Bott-Chern class $c\in H^{1,1}_{BC}(X,\R)$ will be called 
\begin{enumerate}
\item pseudo-effective, if it is
represented by a positive current. 
\item effective, if it decomposes as a finite linear combination $c=\sum_{i=1}^k t_i [D_i]$,
where $t_i\in\R_{\geq 0}$ and $[D_i]$ are Bott-Chern classes associated with  irreducible effective
divisors $D_i$.
 \end{enumerate}    
Similarly, one introduces the notions of  pseudo-effectiveness and effectiveness for a de Rham
cohomology class $c\in H^{1,1}(X,\R)$. 
\end{dt}

We denote by ${\cal P}(X)$,  ${\cal E}(X)\subset H^{1,1}_{BC}(X,\R)$ the cones of
(pseudo-)effective Bott-Chern classes.
\begin{co}\label{pcone} Let $X$ be a complex surface with $b_1(X)$ odd. Then 
\begin{enumerate}
\item  ${\cal P}(X)=C(\{{\gamma_0}\}\cup {\cal D}ou(X)_-^{\rm irr})$ 
\item The inclusion ${\cal E}(X)\subset {\cal P}(X)$ is an equality if and only if $X$ has an
effective divisor representing the trivial real  homology class (i.e. an effective divisor $D$
with $D^2=0$). 
\item If $c\in {\cal P}(X)$, then the de Rham cohomology class of $c$ is effective.
\end{enumerate}
\end{co}
\pf The first statement follows directly from the previous corollary. For the second, note that
when $X$ does not admit any effective divisor representing the trivial real homology class, then
$\gamma_0\not\in {\cal E}(X)$. Indeed, suppose that $\gamma_0$   decomposes in Bott-Chern
cohomology as
$\gamma_0=\sum_{i=1}^k t_i [D_i]$, where $D_i$ are irreducible effective divisors and $t_i>0$. 
Consider the subspaces 
$$L:=\left\{(a_1,\dots a_k)\in\R^k|\ \sum_i a_i [D_i] =0\hbox{ in }
H^{2}(X,\R)\right\}\subset
\R^k
$$
$$L_\Q:=\left\{(a_1,\dots a_k)\in\Q^k|\ \sum_i a_i [D_i] =0\hbox{ in }
H^{2}(X,\R)\right\}\subset
\Q^k
$$
Since the cohomology classes $[D_i]\in H^{2}(X,\R)$ are rational, it follows that
$L=L_\Q\otimes_\Q\R$ so that $L_\Q\cap \Q_{>0}^k$ is dense in $L\cap \R_{>0}^k$. Therefore, we can
find positive rationals
$q_i$ such that $\sum_i q_i [D_i]=0$ in $H^{1,1}(X,\R)$; this obviously gives an   effective
divisor representing the trivial real homology class. The third statement follows from the first.
\qed

\begin{re} There exist  complex surfaces with $b_1(X)$ odd admitting effective
divisors representing the trivial real homology class, but admitting no irreducible effective
divisors with this property.
\end{re}

Indeed, let $X$ be an {\it exceptional compactification of an affine line bundle
over an elliptic curve} \cite{Na1}. Such a surface belongs to class VII and contains
   a cycle
$C=\sum D_i$ of
$b_2(X)$ rational curves $D_i$ having $D_i^2=-2$, $C^2=0$. Every
homological trivial effective divisor of such a surface is a positive integer
multiple of the cycle $C$.
\\ \\

The third statement in Corollary \ref{pcone} above has the following important consequence, which
can be regarded as a  strong   existence criterion for curves on non-K\"ahlerian surfaces.

\begin{co}\label{curve-ex} Let ${\cal L}$ be a holomorphic line bundle over a
complex surface
$X$ with $b_1(X)$ odd. Suppose that $\deg_g({\cal L})\geq 0$ for every Gauduchon
metric
$g$ on $X$. There exists $n\in\N^*$ such that the de Rham Chern class $n
c_1^{DR}({\cal L})\in H^{1,1}(X,\R)$ is represented by an effective divisor.
 \end{co}
\pf  If $\deg_g({\cal L})\geq 0$ for every Gauduchon metric $g$ on $X$, then the Chern class
$c_1^{BC}({\cal L})$ in Bott-Chern cohomology is pseudo-effective, so it decomposes as
$$c_1^{BC}({\cal L})=t_0[\gamma_0]+ \sum_{D\in{\cal D}ou(X)_-^{\rm irr}}   t_D [D]\ .
$$
with coefficients $t_0$, $t_D\geq 0$. Therefore, for the de Rham Chern class, one gets 
\begin{equation}\label{dec}
c_1^{DR}({\cal L})=\sum\limits_{D\in{\cal
D}ou(X)_-^{\rm irr}}   t_D [D]\ .
\end{equation}
On the other hand  $c_1^{DR}({\cal L})$, $c_1^{DR}({\cal O}(D))$, $D\in {\cal
D}ou(X)_-^{\rm irr}$ belong  to the $\Q$-vector space $H^{2}(X,\Q)$. Putting
$d:=\#({\cal D}ou(X)_-^{\rm irr})$, we see that the set $A$ of real systems $(t_D)_{D\in {\cal
D}ou(X)_-^{\rm irr}}$ satisfying (\ref{dec}) is an affine subspace of
$\R^d$ defined by a   linear  system with rational coefficients.  Therefore $A\cap \Q_{\geq 0}^d$
is dense in $A\cap \R_{\geq 0}^d$, so one can find rational non-negative coefficients satisfying
(\ref{dec}).
\qed

\section{Applications}

\subsection{The total Ricci scalar   of a non-K\"ahlerian surface}

Let $X$ be a complex surface. Let $g$ be a Hermitian metric $g$ on $X$, $A_g$
the corresponding Chern connection  on the holomorphic tangent bundle $\Theta_X$, and $s_g$ the 
Ricci scalar   of $g$ which is defined by 
$$s_g:=i\Lambda_g \tr F_{A_g}  
$$
 (see \cite{G}). The total  Ricci scalar   of $g$ is
$$\sigma_g:=\int_X s_g vol_g=\int_X i\omega_g\wedge \tr F_{A_g}\ .
$$
For a Gauduchon metric $g$, one has the following important interpretation of the total Ricci
scalar  
\begin{equation}\label{deg}
\sigma_g=2\pi\deg_g (\Theta_X)=-2\pi\deg_g({\cal K}_X)\ ,
\end{equation}
where  $\Theta_X$ is the holomorphic tangent bundle of $X$, ${\cal K}_X=\det(\Theta_X)^\vee$ is
the canonical line bundle and $\deg_g:\Pic(X)\to\R$ is the Gauduchon degree associated with $g$
(\cite{G}, \cite{LT}).

The purpose of this section is to describe explicitly the set
$$\sigma(X):=\{\gamma_g|\ g\in{\cal G}(X),\ \int_X vol_g=1\}=\{-2\pi\deg_g({\cal K}_X)|\ g\in{\cal
G}(X),\ \int_X \omega_g^2=2\}\ ,
$$
and to decide whether it is a deformation invariant or not.
\begin{lm}\label{sigma} For a Bott-Chern cohomology class $u\in H^{1,1}_{BC}(X)$ on  surface with odd first Betti number
put
$$\sigma(u):=\{ h\cdot u|\ h\in{\cal H}^+(X),\ h^2=1\}\ .
$$
Then 
$$\sigma(u)=\left\{
\begin{array}{ccc}
0&\rm   when&  u=0\\
 (0,\infty)&\rm  when&u\in {\cal P}(X)\setminus\{0\}\\
 (- \infty,0)&\rm  when&u\in -{\cal P}(X)\setminus\{0\}\\
 (-\infty,\infty)&\rm when &u\not\in{\cal P}(X)\cup (-{\cal P}(X))\ .
\end{array}\right.
$$
\end{lm} 
 \pf An element $u\in {\cal P}(X)\setminus\{0\}$  
obviously satisfies $\sigma(u)\subset (0,\infty)$.  For the converse inclusion we proceed as follows:

Let $\omega_0$ be the K\"ahler form of a fixed Gauduchon
metric $g_0$. Our description of the positive cone ${\cal H}_+$ shows that the whole half-line
$$\left\{[\omega_0]+t\gamma_0\left|\ ([\omega_0]+ 
t\gamma_0)^2>0\right\}\right.=\left\{[\omega_0]+t\gamma_0\left|\  [\omega_0]^2+ 
2t\omega_0\cdot\gamma_0>0\right\}\right.
$$
is contained in ${\cal H}_+$. Put $h_t:=[\omega_0]+t\gamma_0$ for $t>
-\frac{\omega_0^2}{2\omega_0\cdot\gamma_0}$. It suffices to notice that, for any $u\in{\cal
P}(X)\setminus\{0\}$, one has
$$\left\{\frac{1}{\sqrt{h_t^2}}\ h_t\cdot u\ \left|\  
t> -\frac{\omega_0^2}{2\omega_0\cdot\gamma_0}  \right.\right\}=$$
$$=\left\{\frac{1}{\sqrt{[\omega_0]^2+ 
2t\omega_0\cdot\gamma_0}}\ \omega_0\cdot u\
\left|\ t>
 -\frac{\omega_0^2}{2\omega_0\cdot\gamma_0} \right. 
\right\}=(0,\infty)\ . 
$$
This proves the first three equalities. 

For the fourth, suppose that $H^{1,1}_{BC}(X,\R)\ni u\not\in {\cal P}(X)\cup(-{\cal
P}(X))$. By Corollary \ref{equiv} there exist   Gauduchon metrics $g_1$, $g_2$ such
that the corresponding K\"ahler forms satisfy
$$[\omega_1]\cdot u>0\ ,\ [\omega_2]\cdot u<0
$$
Modifying the  two classes $[\omega_1]$, $[\omega_2]$ by $t\gamma_0$ as   above, one gets easily two half-lines $l_1$,
$l_2\subset {\cal H}_+$ such that $l_1\cdot u=(0,\infty)$,  $l_2\cdot u=(-\infty,0)$.
\qed 

By Lemma \ref{sigma} and formula (\ref{deg}), the set $\sigma(X)$ is determined by the position of the Chern class
$c_1^{BC}({\cal K}_X)$ in Bott-Chern cohomology with respect to the cones $\pm{\cal P}(X)=\pm C(\{\gamma_0\}\cup {\cal
D}ou(X)^{\rm irr}_-)$.
\\

Taking into account that     
 algebraic surfaces of Kodaira dimension $-\infty$ allow K\"ahler metrics with positive total
scalar curvature, the following remark might be surprising:

\begin{re} Let $X$ be any class $VII$ surface whose minimal model is an Inoue surface. Then
the class
$c_1^{BC}({\cal K}_X)$ is non-zero and pseudo-effective, hence 
$\sigma(X)=(-\infty,0)$.
\end{re}
\pf  Let $\H$ be the upper half-plane. An Inoue  surface is a quotient of $\H\times\C$ by  a properly
discontinuous group
$G$ of affine transformations. There are three classes of Inoue surfaces \cite{In}, denoted by $S_M^\pm$,
$S^+_{N,p,q,r,t}$,
$S^-_{P,p,q,r}$. Here $M\in SL(3,\Z)$ is a matrix with a single  real eigenvalue $\alpha>1$, $N\in SL(2,\Z)$
has two positive real eigenvalues 
$\alpha^{-1}, \alpha>1$,   $P\in GL(2,\Z)$  has two real eigenvalues $\alpha>1$ and $-\alpha^{-1}$. The
  symbols $r$, $t$ denote numbers
$r\in \Z\setminus\{0\}$, $t\in\C$ whereas $p$, $q$ are real numbers  satisfying a certain
integrality condition.

Taking into account the way in which the group acts on pairs $(w,z)\in\H\times\C$, one checks easily that in the case of the surfaces $S^+_{N,p,q,r,t}$, the form $\frac{dw\wedge dz}{{\rm Im}(w)}$ descends to a differentiable nowhere vanishing $(2,0)$-form  on  $S$. This shows that  setting 
$$h_{(w,z)}(dw\wedge dz)={\rm Im}(w)^2$$
 one gets a Hermitian metric on  the line bundle ${\cal K}_S$. Therefore,  setting $w=u+iv$,
we see that the form 
$$\frac{i}{\pi}\bar\partial \ \frac{\partial v}{v}=\frac{i}{\pi}\bar\partial(-\frac{i}{v}
dw)=\frac{1}{\pi}\bar\partial(\frac{1}{v} dw)=\frac{i}{\pi}(-\frac{1}{v^2}d\bar w\wedge dw)=\frac{i}{\pi
v^2} dw\wedge d\bar w
$$
descends to a closed $(1,1)$-form representing the  Chern class of ${\cal K}_S$ in  Bott-Chern
cohomology. But this form is a non-trivial positive current. This shows that $c_1^{BC}({\cal K}_S)$
is non-zero and pseudo-effective.  For an Inoue surface $S$ of type 
$S^-_{P,p,q,r}$ the formula $h_{(w,z)}(dw\wedge dz)={\rm Im}(w)^2$ still defines  a Hermitian metric on ${\cal K}_{S}$ whose Chern form is positive. For $S_M^\pm$ on takes   $h_{(w,z)}(dw\wedge dz)={\rm Im}(w)$\footnote{ I am indebted to V. Apostolov and  G. Dloussky for pointing out that the surfaces  $S^\pm_M$, $S^-_{P,p,q,r}$ require a slightly different argument}.  

For a blown up Inoue surface $X\textmap{p} S$, one just notices that 
$$c_1^{BC}({\cal
K}_X)=p^*( c_1^{BC}({\cal K}_S))+ [E]\ ,$$
 where $E$ is an effective divisor.
\qed

For a Hopf surface, one has: 
\begin{re} Any primary Hopf $H$ has an anti-canonical divisor. Therefore for  such a surface
one has
$\sigma(H)=(0,\infty)$.
\end{re}
\pf A primary a Hopf surface $H$ of the form   $\C^2\setminus\{0\}/\langle
T\rangle$, where
$$T:(z_1,z_2)\mapsto (\alpha_1 z_1,\alpha_2 z_2)$$
 (where $0<|\alpha_1|\leq |\alpha_2|<1$)
has ${\cal K}_H={\cal O}_H(-C_1-C_2)$, where $C_i$ are the elliptic curves defined by the
equations $z_i=0$. If $T$ has the form 
$$(z_1,z_2)\mapsto (\alpha_1 z_1+az_2^m,\alpha_2 z_2)$$
where $\alpha_2^m=\alpha_1$, one has ${\cal K}_H={\cal O}_H(-(m+1)C)$ where $C$ is the
elliptic curve defined by the equation $z_2=0$.
\qed

For a blown up Hopf surface the result is more complicated.
\begin{pr} Let $X$ be class $VII$ surface with $b_2(X)>0$ whose minimal model is a primary
Hopf surface. Then $c_1^{BC}({\cal K}_X)\not\in {\cal P}(X)\cup(-{\cal P}(X))$. In
particular $\sigma(X)=(-\infty,\infty)$.
\end{pr}
\pf  For simplicity we give the proof only for a single blow up. Let $H$ be a primary
Hopf surface with anti-canonical effective divisor $A$, let $\pi:X\to H$ the
blow  up at a point $x_0\in H$ with exceptional divisor $E$, and denote by $\tilde A$ the
proper transform of $A$. Then
${\cal K}_X$ decomposes as
$${\cal K}_X={\cal O}_X(-D)\otimes {\cal O}_X(E)\ ,
$$
where   ${\cal O}(D)=\pi^*({\cal O}(A))$, $D=\tilde A+  k E$ (where $k\geq 0$ is the
incidence order between $x_0$ and $A$).   

Since $D$ is homologically trivial, we get
$$c_1^{BC}({\cal K}_X) = -t_0\gamma_0+  [E] \ .
$$
Since $E$ is the only irreducible effective divisor with negative self-intersection and 
$\gamma_0$, $[E]$ are linearly independent in $H^{1,1}_{B,C}(X,\R)$, we see easily that
$c_1^{BC}({\cal K}_X)\not\in {\cal P}(X)\cup(-{\cal P}(X))$.
\qed
\\
{\bf Remark:} There is standard way to endow a blown up surface $\pi:\hat X_{x_0}\to X$ with
a Gauduchon metric (see \cite{Bu1},  \cite{LT}). The idea is to choose a Gauduchon
metric
$g$ on $X$ and to note that there exists a closed $(1,1)$-form $\eta$ representing the class
of the exceptional curve $E$ whose restriction to this curve is the opposite of its
Fubiny-Study volume form. It will follow that, for all sufficiently small $\varepsilon>0$,
the form $\pi^*(\omega_g)-\varepsilon\eta$ is positive and $i\bar\partial\partial$-closed,
so it corresponds to a Gauduchon metric $\hat g_\varepsilon$ on $\hat X_{x_0}$. 

The volume of the exceptional divisor  with respect to a metric $\hat
g_\varepsilon$  is small. Therefore, in this way one gets Gauduchon metrics with {\it
positive} total Ricci scalars on blown up Hopf surfaces; it is not clear at all how to
construct explicitly Gauduchon metrics with {\it negative} total Ricci scalars on these
surfaces.
\\ \\
For the minimal case one has:

\begin{pr}\label{nonpse} Let $X$ be a minimal class VII surface with $b_2(X)>0$.
The  class
$c_1^{BC}({\cal K}_X)$ cannot be pseudo-effective. Therefore, such a surface has
either
$\sigma(X)=(-\infty,\infty)$ (when $c_1^{BC}(\Theta_X)$ is not pseudo-effective) or
$\sigma(X)=(0,\infty)$ (when $c_1^{BC}(\Theta_X)$ is   pseudo-effective).
\end{pr}
\pf Suppose that  $c_1^{BC}({\cal K}_X)$ was pseudo-effective. By Corollary
\ref{curve-ex}, there exists $n\in\N^*$ such that $nc_1^{DR}({\cal K}_X)=PD([E])$
for an effective divisor $E\subset X$. This gives $\langle c_1^{DR}({\cal
K}_X),[E]\rangle=\frac{1}{n} c_1^{DR}({\cal
K}_X)^2=-\frac{1}{n} b_2(X)<0$, which contradicts Lemma 1.1 in \cite{Na3}.
\qed

There exist many minimal class  VII surfaces with
pseudo-effective  $c_1^{BC}(\Theta_X)$, for instance the surfaces 
allowing a pluri-anticanonical divisor. A hyperbolic Inoue surface
$X$
\cite{Na1}, \cite{Na2} has two cycles
$A$, $B$ of rational curves, and one has ${\cal K}_X\simeq{\cal O}_X(-A-B)$  (see
Lemma 2.8 in
\cite{Na1} and the proof of Lemma 4.7). Similarly, a  {\it half Inoue surface}
$X$
\cite{Na1} has a cycle $C$ of
$b_2(X)$ rational curves and an order two flat line bundle ${\cal L}$ such that ${\cal
K}_X\simeq {\cal L}\otimes{\cal O}_X(-C)$; thus $2C$ is a bi-anti-canonical
divisor.\\

There  also exist minimal class  $VII$-surfaces  $X$  with
\begin{enumerate}
\item   pseudo-effective  
Bott-Chern class Chern class $c_1^{BC}(\Theta_X)$ but allowing no
pluri-anticanonical divisors, 
\item   non-pseudo-effective 
$c_1^{BC}(\Theta_X)$.
\end{enumerate}

Any {\it known} minimal class VII surface with $b_2>0$ is   the special fibre
$X_0$ of a family of surfaces ${\cal X}\to D$  whose fibres $X_t$, $t\ne 0$ are all blown up
primary Hopf surfaces. If the GSS conjecture was true (which has been proved for  $b_2=1$
\cite{Te2}),
  {\it any} minimal class VII surface with $b_2>0$ would be a degeneration of a family of
blown up primary Hopf surfaces.
\begin{co}The assignment $X\mapsto \sigma(X)\subset\R$ is not a deformation invariant for
class $VII$ surfaces. More precisely there exist families ${\cal X}\to D$ of such surfaces
such that $\sigma(X_t)=(-\infty,\infty)$ for any $t\ne 0$ and $\sigma(X_0)=(0,\infty)$.
\end{co}
\pf It suffices to consider a one parameter family of blown up primary Hopf surfaces
degenerating to a minimal class VII surface with   pseudo-effective
$c_1^{BC}(\Theta_X)$ (for instance a hyperbolic Inoue surface).

\subsection{The stability of the canonical extension of a class
$VII$-surface}

Class VII surfaces are not completely classified yet. The main obstacle is the
``Global Spherical Shell (GSS) conjecture"   (\cite{Na2}, p. 220) which states that
any minimal class VII surface $X$ with $b_2(X)>0$ contains a global spherical
shell, i.e. an open submanifold $S\subset X$ biholomorphic to a neighborhood of
$S^3$ in $\C^2$  such that $X\setminus S$ is connected.  Minimal class VII surfaces
containing a global spherical shell  are well understood; any such surfaces $X$
contains
$b_2(X)$ rational curves, but there are many possibilities for the intersection
graph of these curves.  This intersection graph is {\it not} a deformation
invariant.  By a fundamental result of Dloussky-Oeljeklaus-Toma
\cite{DOT}, any minimal class VII surface $X$ which has $b_2(X)$ rational curves
does contain a GSS, so the classification of class VII surfaces reduces to the  
question: `` does any minimal class $VII$ surface with
$b_2(X)>0$ possess $b_2(X)$ rational curves''?
\\ \\ 

Let $X$ be an arbitrary  class $VII$ surface. By Serre duality    $h^1({\cal
K}_X)=1$, so there exists a (up to isomorphy) unique    rank 2-holomorphic bundle
${\cal A}$ which is the central term of a nontrivial extension
$$0\map {\cal K}_X\textmap{s} {\cal A}\textmap{t} {\cal O}_X\map  0\ ,
$$
which will be called  {\it the canonical extension of} $X$. The problem treated
in this section is: does there exist  Gauduchon metrics on $X$ with respect
to which
${\cal A}$ is stable?  The problem is not easy: when $\deg_({\cal K}_X)<0$, the
obvious  line subbundle  ${\cal K}_X$ of
${\cal A}$ does not destabilize it, but it is of course possible that ${\cal A}$ is
destabilized by  another line bundle. This would imply that ${\cal A}$ can be
written as extension in a different way. On the other hand, we will see that {\it
writing   a rank 2-bundle as an extension in two different ways, implies the
existence of effective divisors with very special properties.} Therefore, the
stability of ${\cal A}$ is related to the existence of curves on the base
manifold $X$.   This is an important remark  because, by  
Dloussky-Oeljeklaus-Toma's theorem,  the GSS conjecture   reduces to
the existence of ``sufficiently many" curves.
\\ \\
{\bf Example:} Let $X$ be an Inoue-Hirzebruch surface (a hyperbolic Inoue surface)
\cite{Na1}. Such a surface has two disjoint cycles $A$, $B$ of rational curves, and
${\cal K}_X\simeq {\cal O}_X(-A-B)$ (\cite{Na1} p. 402, 419).
We state that, in this case one has ${\cal A}\simeq {\cal O}(-A)\oplus {\cal O}(-B)$, so the
canonical extension of such a surface is non-stable with respect to any Gauduchon
metric. Indeed, one has ${\cal K}_X^\vee\otimes [{\cal O}(-A)\oplus {\cal
O}(-B)]={\cal O}(B)\oplus{\cal O}(A)$  so, since $A\cap B=\emptyset$, one gets a
bundle embedding ${\cal K}_X\hookrightarrow {\cal O}(-A)\oplus {\cal O}(-B)$.
Therefore, ${\cal O}(-A)\oplus {\cal O}(-B)$ is an extension of ${\cal O}_X$ by
${\cal K}_X$, and this extension cannot be trivial because $H^0({\cal O}(-A)\oplus
{\cal O}(-B))=0$.
\\

Let $(e_1,\dots e_{b_2(X)})$ be a basis of $ {H^2(X,\Z)}/{\rm Tors}$ such that $e_i^2=-1$
and $c_1^\Q({\cal K}_X)=\sum_ i e_i$. the existence of such a basis follows easily
(see \cite{Te2}) from Donaldson's   theorem on smooth manifolds with definite
intersection form
\cite{Do}. For a subset $I\subset\{1,\dots,b_2(X)\}$ we put  
$$e_I:=\sum_{i\in I}
e_i\ ,\ \bar I:=\{1,\dots,b_2(X)\}\setminus I\ .$$
\begin{lm}\label{topo} Let ${\cal E}$ be any holomorphic 2-bundle with $\det({\cal E})={\cal
K}_X$,
$c_2({\cal E})=0$ and let $j:{\cal L}\hookrightarrow {\cal E}$ a rank 1 locally free subsheaf
with torsion free quotient. Then $j$ is a bundle embedding (i.e. fibrewise injective) and
there exists a subset
$I\subset
\{1,\dots,b_2(X)\}$ such that
$c_1^\Q({\cal L})=e_I$.
\end{lm}
\pf 
The inclusion ${\cal L}\hookrightarrow {\cal E}$  fits in an exact sequence
$$0\map{\cal L}\textmap{j}{\cal E}\textmap{k}{\cal K}_X\otimes {\cal L}^{-1}\otimes {\cal
I}_Z\map 0\ ,
$$
where $Z\subset X$ is a codimension 2 locally complete intersection. Decomposing
$c_1^\Q({\cal L})=\sum _i a_i e_i$ (with $a_i\in\Z$), this gives
$$0=c_2({\cal E})=|Z]+\sum a_i(a_i-1)\ ,
$$
which happens iff $Z=\emptyset$  and $a_i\in\{0,1\}$ for all $i\in\{1,\dots,b_2(X)\}$.
\qed

\begin{pr}\label{teh} Let $S$ be an arbitrary complex surface, let
\begin{equation}
\label{ext}0\map {\cal L}\textmap{a} {\cal E}\textmap{b} {\cal O}_S\map 0
\end{equation}
an exact sequence, and $\varepsilon:=\delta_h(1)\in H^1({\cal L})={\rm Ext}^1({\cal
O}_X,{\cal L})$ the corresponding extension invariant, where $\delta_h$ stands for the
connecting operator in the associated cohomology sequence. Let 
$D\subset X$ a (possibly empty, possible non-reduced) effective divisor, and $u:{\cal
O}_S(-D)\to {\cal O}_S$ the canonical morphism. Let
$V\subset
\Hom({\cal O}_S(-D),{\cal E})=H^0({\cal E}(D))$   be the set of liftings of $u$ to ${\cal
E}$. Then 

\begin{enumerate}

\item If non-empty, $V$ is an affine space modeled over the vector space $H^0({\cal L}(D))$.

\item The restriction $\resto{v}{D}\in H^0({\cal E}_D(D))$ of any lifting $v\in V$ to $D$
belongs to  the subspace $H^0({\cal L}_D(D))$  of $H^0({\cal E}_D(D))$, so it defines an
element 
$$\rho(v)\in  H^0({\cal L}_D(D))\ .$$
\item For every $v\in V$, the element $\rho(v)$ is a lifting of $\varepsilon$ via the
connecting operator
$$\delta_v:H^0({\cal L}_D(D)) \map H^1({\cal L})
$$
associated with the exact sequence $0\to{\cal L}\to{\cal L}(D)\to{\cal L}_D(D)\to 0$.
\item The map $V\ni v\mapsto \rho(v)\in H^0({\cal L}_D(D))$ defines a bijection between the
quotient
$V/H^0({\cal L})$ and the space $H_\varepsilon$ of $\delta_v$-liftings of $\varepsilon$ in 
$H^0({\cal L}_D(D))$. Here   $H^0({\cal L})$ was regarded as a   subspace of the model
vector space
$H^0({\cal L}(D))$.
\item The vanishing locus $Z(v)\subset S$ of $v$ is contained in $D$ and coincides with the
vanishing locus $Z(\rho(v)))\subset D$. In particular, the lifting
$v:_S{\cal O}(-D)\to{\cal E}$ of $u$ is a bundle embedding if and only if
$\rho(v)$ is a trivialization of the line bundle ${\cal L}_D(D)$ over $D$.
\end{enumerate}
\end{pr}
\pf  The first statement is obvious. For the second, use the following  sheaf
diagram with exact rows and exact columns

\begin{equation}\label{diag}
\begin{array}{ccccccccc}
&&0\phantom{iii}&&0\phantom{iii}&&0\phantom{iii}&&\\
\\
&&\downarrow\phantom{ii }&&\downarrow\phantom{ii }&&\downarrow\phantom{ii i}&&\\
\\
0&\map& {\cal L}\phantom{ii} &\textmap{a} &{\cal E}\phantom{ii} &\textmap{b} &{\cal
O}_S\phantom{ii}&\map& 0\\
\\
 &&\downarrow i&&\downarrow i'&&\downarrow i''&&\\
\\
0&\map& {\cal L} (D)&\textmap{a'} &{\cal E} (D)&\textmap{b'} &{\cal O}_S(D)&\map& 0\\
\\
&&\downarrow p&&\downarrow p'&&\downarrow p''&&\\
\\
0&\map& {\cal L}_D(D)&\textmap{a''} &{\cal E}_D(D)&\textmap{b''} &{\cal O}_D(D)&\map& 0\\
\\
&&\downarrow\phantom{ii}&&\downarrow\phantom{ii}&&\downarrow\phantom{ii}&&\\
\\
&&0\phantom{iii}&&0\phantom{iii}&&0\phantom{iii}&&
\end{array}
\end{equation} 
and note that the image of $u\in H^0({\cal O}_S(D))$ in $H^0({\cal O}_D(D))$ vanishes, so
the image of $v$ in ${\cal E}_D(D)$ belongs to ${\cal L}_D(D)$.\\
\\ \\
3. The first row and the first column in  the diagram (\ref{diag}) can be regarded
as
 resolutions of the rank 1 locally free sheaf ${\cal L}$. This diagram also yields a third
resolution of the same sheaf, namely the simple (or total) complex associated with
the
 double complex (\ref{diag}).
$$0\to {\cal L}\textmap{(a,i)}  {\cal E}\oplus {\cal L}(D)\stackrel{A}{\to} {\cal
O}_S\oplus{\cal E}(D)\oplus {\cal L}_D(D)\stackrel{B}{\to} {\cal O}_S(D)\oplus{\cal
E}_D(D)\stackrel{C}{\to} {\cal O}_D(D)\to 0\ .
$$
Truncating this resolution, one gets the short exact sequence
\begin{equation}\label{us}
0\to {\cal L}\textmap{(a,i)}  {\cal E}\oplus {\cal
L}(D)\textmap{A}  
\im(A)=\ker(B)\to 0\ . 
\end{equation}
The idea of the proof is to notice that a lift $v\in V$ of $u$ defines an element $r(v)$ in
$H^0(\im(B))$, namely
$$r(v)=(1,v,\rho(v))\ .
$$
Let $\partial$ be the connecting operator associated with the short exact sequence
(\ref{us}). One can compute $\partial(r(v))\in H^1({\cal L})$ in two ways: comparing the
exact sequence (\ref{us}) with the first row in (\ref{diag}) and using the
functoriality of the connecting operator, one gets 
$\partial(r(v))=\delta_h(1)=\varepsilon$, whereas comparing (\ref{us}) with the first column
in (\ref{diag}), one has $\partial(r(v))=\delta_v(v)$.
\\ \\
4. If $\rho(v)=\rho(v')$, the $v-v'\in H^0({\cal E}_D(D))$ is mapped to 0 via both vertical
an horizontal  arrows in (\ref{diag}). A simple diagram chasing shows that $v-v'$ comes from
$H^0({\cal L})$ via the obvious morphism. This proves the injectivity.  Let now $w\in
H^0({\cal L}_D(D))$ be an element which is mapped to $\varepsilon$ via $\delta_v$.\\

For surjectivity, let $w\in H^0({\cal L}_D(D))$ be a lift of $\varepsilon$ via $\delta$.
Since $\varepsilon=\delta_v(w)$, it follows that the image of $\varepsilon$ in $H^1({\cal
L}(D))$ vanishes. Similarly, since $\varepsilon=\delta_h(1)$,   the image of
$\varepsilon$ in $H^1({\cal E})$ will vanish, too.  Therefore, in the cohomology sequence 
associated with (\ref{us}),  
$\varepsilon$ is mapped to 0 in ${\cal E}\oplus {\cal
L}(D)$, so $\varepsilon$ can be written as $\partial(\theta)$, for an element
$\theta=(\varphi,\psi,\chi)\in H^0(\ker(B))$. Using again the functoriality of the
connecting operator, we see that
$\delta_h(\varphi)=\delta_v(\chi)=\varepsilon=\delta_h(1)=\delta_v(w)$. We can modify the
triple
$\theta$ by a suitable element in $A(H^0({\cal E})\oplus H^0({\cal L}(D)))$  to
get  a lift $\theta'$ of $\varepsilon$ via $\partial$, having the first component 1 and
the third component $w$. The second component $v$ of $\theta'$ will satisfy $\rho(v)=w$.\\
\\
5. Consider, in general,  an epimorphism   $\pi:{\cal F}\to  {\cal G}$ of holomorphic
vector bundles over a complex space $Y$, and let $\sigma$ be a holomorphic section of ${\cal
F}$. the vanishing locus $Z(\sigma)$ of $\sigma$ is the complex subspace of $Y$ defined by
the ideal sheaf $\sigma^\vee({\cal F}^\vee)\subset {\cal O}_Y$. One has the following useful 
associativity property:
$$Z(\sigma)=Z(\resto{\sigma}{Z(\pi\circ \sigma)})\ , 
$$
where the restriction $\resto{\sigma}{Z(\pi\circ \sigma)}$ can be regarded as a section in
the holomorphic bundle $\resto{\ker(\pi)}{Z(\pi\circ \sigma)}$.

The result follows by applying this associativity principle to the epimorphism $b':{\cal
E}(D)\to {\cal O}_S(D)$ and to notice that $Z(b'\circ v)=Z(u)=D$. For the second statement,
note that ${\cal L}_D(D)$ is a line bundle on $D$, so the vanishing locus of a section in
${\cal L}_D(D)$ is empty if and only if defines a global trivialization of this line bundle
(the condition that it does not  vanish at any point is not sufficient for a non-reduced
divisor $D$).
\qed
\begin{co} With the notations and in the conditions of Proposition \ref{teh}, the
natural map
${\cal O}(-D)\to {\cal O}_X$ can be lifted to a bundle embedding ${\cal
O}(-D)\hookrightarrow {\cal E}$   if and only if
there exists a section $\alpha\in H^0({\cal L}\otimes{\cal O}_D(D))$ with the
properties
\begin{enumerate}
\item $\alpha$ defines a trivialization of ${\cal L}\otimes{\cal O}_D(D)$,
\item The image of $\alpha$ in $H^1({\cal L} )$ via the connecting operator $\delta_v$ is the
invariant
$\varepsilon$ of the given extension (\ref{ext}).
\end{enumerate}
\end{co}
\begin{co}\label{secext} Let $X$ be a minimal class VII surface with $b_2(X)>0$ and ${\cal
A}$ its canonical extension. The bundle
${\cal A}$ can be written as an extension
$$0\map {\cal M}\map {\cal A}\map {\cal K}_X\otimes {\cal M}^{-1}\map 0
$$
iff and only if there exists a {\it non-empty} effective divisor $D\subset  X$ satisfying the
following properties:
\begin{enumerate}
\item ${\cal M}\simeq {\cal O}_X(-D)$,
\item ${\cal K}_X\otimes {\cal O}_D(D)\simeq {\cal O}_D$.
\item $h^0({\cal K}_X\otimes {\cal O}_D(D))-h^0({\cal K}_X\otimes {\cal O}(D))=1$
\end{enumerate}
\end{co}
\pf

If  ${\cal M}$ is a line subbundle of ${\cal A}$, it must admit a non-trivial map to ${\cal
O}_X$, so it must be isomorphic to ${\cal O}(-D)$ for an effective divisor $D\subset X$.
$D$ must be non-empty, because the extension defining ${\cal A}$ does not split. The
second   condition is necessary by Corollary
\ref{secext}; we show that the third condition is also necessary. In our case, the cohomology
exact sequence associated with the first vertical column in (\ref{diag}) reads
$$H^0({\cal K}_X)=0\map H^0({\cal K}_X\otimes{\cal O}_X(D))\map H^0({\cal
K}_X\otimes {\cal O}_D(D))\textmap{\delta_v} H^1({\cal K}_X)\simeq\C
$$

 On the other hand, by Corollary \ref{secext}, the  map $\delta_v$ must be
non-trivial, so one has 
$h^0({\cal K}_X\otimes {\cal O}_D(D))-h^0({\cal K}_X\otimes {\cal O}(D))=1$.
\\

In order to prove that the three conditions are also sufficient,    it suffices to
show that they imply the existence of  a trivialization of  ${\cal K}_X\otimes
{\cal O}_D(D)$ which is mapped onto a non-trivial element of the line $H^1({\cal
K}_X)$. 

Since $a(X)=0$ (\cite{Na3} p. 477), one has  $h^0({\cal K}_X\otimes {\cal
O}(D))\leq 1$  so, by the third condition one has either $h^0({\cal K}_X\otimes
{\cal O}_D(D))=1$ and $h^0({\cal K}_X\otimes {\cal O}(D))=0$ or $h^0({\cal
K}_X\otimes {\cal O}_D(D))=2$ and $h^0({\cal K}_X\otimes {\cal O}(D))=1$. Taking
into account that
${\cal K}_X\otimes {\cal O}_D(D)\simeq {\cal O}_D$, the claim is obvious in the
first case  because, in this case, {\it any} non-trivial
section of ${\cal K}_X\otimes {\cal O}_D(D)$ will be a trivialization  which is
mapped onto  a non-trivial element of $H^1({\cal K}_X)$.

In the case  $h^0({\cal K}_X\otimes {\cal O}_D(D))=2$, $h^0({\cal K}_X\otimes {\cal
O}(D))=1$, it suffices to notice that both $\ker(\delta_v)$ and the subset $F$  of $H^0({\cal
K}_X\otimes {\cal O}_D(D))$ consisting of sections which are not trivializations are  
proper Zarisky   closed subsets of the vector space $H^0({\cal
K}_X\otimes {\cal O}_D(D))$, so the complement of their union is non-empty. 

\qed

Let $X$ be a minimal class VII surface. We agree to call {\it a cycle} of $X$ any
reduced divisor $C$ which is either an elliptic curve, or a singular rational
surface with a node, or a cycle of smooth rational curves. In the last two cases
$C$ will be called a cycle of rational curves. In all three cases  the canonical
sheaf
$\omega_C:={\cal K}_X\otimes{\cal O}_C(C)$ is trivial.

\begin{pr}\label{destlb}  Let $X$ be a minimal class VII surface with $b_2(X)>0$ and
let
$D\subset X$ be a non-empty effective divisor of $X$ satisfying
\begin{enumerate}
\item   The canonical line bundle $\omega_D:={\cal K}_X\otimes
{\cal O}_D(D)$ is trivial.
\item $h^0({\cal K}_X\otimes {\cal O}_D(D))- h^0({\cal K}_X\otimes {\cal O}(D))=1$. 
\end{enumerate} 
Then there exists $I\subset\{1,\dots,b_2(X)\}$ such that $c_1^\Q({\cal O}(-D))=
e_I$ and one of the following holds
\begin{enumerate}
\item $D$ is a cycle,
\item ${\cal O}(-D)={\cal K}_X$ (i.e. $D$ is an anti-canonical divisor).  
 
\end{enumerate}
\end{pr}
\pf  The existence of $I\subset \{1,\dots,b_2(X)\}$ such that $c_1^\Q({\cal
O}(-D))= e_I$ follows from Lemma \ref{topo} and Corollary \ref{secext}.
\\  

Since $h^0(\omega_D)=h^1({\cal O}_D)\geq 1$, one gets $h^1(D_{\rm red})\geq 1$ by
Lemma 2.7 in \cite{Na1}. Let $0<C\leq D_{\rm red}$ be   minimal with the property
$h^1(C)\geq 1$.  By Lemma  2.3, and Lemma 2.12 in \cite{Na1} $C$ is either a cycle or
a union of two disjoint cycles. Decompose $D$ as
$D=C+E$ for an effective divisor $E$.  Put
${\cal M}:={\cal K}_X\otimes {\cal O}(D)$. Noting that $h^2({\cal M})=h^0({\cal
O}(-D))=0$, we get the exact sequence
$$0\to H^0({\cal M}(-C))\to H^0({\cal M})\to H^0({\cal M}_C)\to 
$$
$$
\to H^1({\cal
M}(-C))\to H^1({\cal M})\to  H^1({\cal M}_C)\to H^2({\cal
M}(-C))\to   0
$$
\\ 
Case 1. $h^0({\cal K}_X\otimes {\cal O}_D(D))=1$ and $ h^0({\cal K}_X\otimes {\cal
O}(D))=0$.\\

In this case  we get by Riemann-Roch theorem  $h^1({\cal M})=0$, hence  (recalling
that ${\cal M}$ is trivial on
$D$, hence also on $C$)
$$h^2({\cal
M}(-C))=h^1({\cal M}_C)=h^1({\cal O}_C)\geq 1\ .$$
  Therefore $h^0({\cal O}(-E))\geq 1$,
which shows that $E$ is empty, so $D=C$.  In this case $D=C$  must be a single
cycle, because otherwise one would have $h^0({\cal K}_X\otimes{\cal O}_D(D))=2$. 
\\ \\
Case 2. $h^0({\cal K}_X\otimes {\cal O}_D(D))=2$ and $h^0({\cal K}_X\otimes {\cal
O}(D))=1$.\\

Note first that $c_1^\Q({\cal K}_X\otimes {\cal
O}(D))=e_{\bar I}$ and $h^0({\cal K}_X\otimes {\cal
O}(D))>0$. Using Lemma 2.3 in \cite{Na3} (which holds for any minimal
class VII surface  with positive $b_2$) one gets $\bar I=\emptyset$, so 
\begin{equation}\label{Kdec} {\cal K}_X\otimes {\cal
O}(D)={\cal O}(F)
\end{equation}
for an  effective divisor $F$ with $F\cdot F=0$. Therefore, $D$ must be a
{\it numerically} anti-canonical divisor. If $X$ contains no homologically trivial
effective divisors, then $F$ must be empty, so   $D$ is anti-canonical as claimed.\\

Suppose  now that $X$ does contain  homologically trivial divisors. 
Minimal class
$VII$-surfaces with $b_2>0$ containing homologically trivial effective divisors are
classified. Any such surface is an exceptional compactification  of an affine line
bundle  over an elliptic curve
\cite{Na1},
\cite{Na2}, \cite{Na3}, contains a GSS, and contains a homologically trivial cycle
$C$ of $b_2(X)$ rational curves $D_i$. An   exceptional
compactification of a non-linear affine line bundle does not contain any other curve
but the  irreducible components $D_i$ of $C$, which all satisfy  $\langle c_1({\cal
K}_X), D_i\rangle =0$. Therefore, on such a surface there exist no anti-canonical
numerically divisor. 

The  exceptional compactifications of linear line bundles are called parabolic Inoue
surfaces. Such a surface  contains a smooth elliptic curve $Z$ with $Z\cdot 
Z=-b_2(X)$,
$Z\cap C=\emptyset$. In this case one has
${\cal K}_X={\cal O}_K(-C-Z)$ and the only homologically trivial effective
divisors are $nC$, $n\in\N$. Therefore (\ref{Kdec}) would imply a linear
equivalence of the form   $D\sim (n+1)C+Z$. Since $a(X)=0$, one has $D= (n+1)C+Z$.

We claim that only for $n=0$ one can have $\omega_D={\cal O}_D$, which   will
complete the proof. Indeed, if
$\omega_{(n+1)C}={\cal O}_{(n+1)C}$  then, taking into account $\omega_C={\cal
O}_C$,  one would have  ${\cal O}_C(nC)={\cal O}_C$.

But ${\cal O}(C)$ is a flat line bundle on $X$ which is associated
with a representation $\rho:\Z\simeq \pi_1(X)\to \C^*$ with $|\rho(1)|\ne 1$ (see
(\cite{D} section 1.2)\footnote {I am indebted to Georges Dloussky, who kindly
explained me this important property of the line bundle
${\cal O}_X(C)$.}. On the other hand,
the natural map
$H_1(C,\Z)\to H_1(X,\Z)$ is an isomorphism (\cite{Na1} p. 404). Therefore, for
$n\geq 1$, the restriction of
${\cal O}(nC)$ to $C$ is a flat line bundle  associated with a nontrivial
representation $\Z\simeq
\pi_1(C)\to \C^*$, so $H^0({\cal O}_C(nC))=0$.
 \qed

Therefore, a line subbundle of ${\cal A}$ is either  isomorphic to ${\cal K}_X$
or to a line bundle of the form ${\cal O}(-C)$ for a cycle $C\subset X$. We can prove
now:

\begin{thry}\label{intth} Let $X$ be a minimal class VII surface with $b_2>0$.
Suppose that
${\cal A}$ is unstable for any Gauduchon metric on
$X$. Then   one of the following holds:
\begin{enumerate}
\item $X$ contains two cycles, i.e. $X$ is either a hyperbolic or a parabolic Inoue
surface \cite{Na1}.

In this case ${\cal A}$ is a direct sum of line bundles.
\item $c_1^{BC}(\Theta_X)$ is pseudo-effective, $X$ contains a single cycle $C$ and
the Bott-Chern class
$c_1^{BC}({\cal K}_X^\vee(-2C))$ is  pseudo-effective.

\end{enumerate}
\end{thry}
\pf Let  $X$ be a minimal class VII surface with $b_2(X)>0$ which does not contain
two cycles. In other words, $X$ is neither a hyperbolic nor a parabolic Inoue
surface. We will prove that, if
${\cal A}$ is unstable with respect to any Gauduchon metric on $X$  then 2. holds. 
By Proposition \ref{nonpse} we have to consider only the following two cases:
\begin{enumerate}
\item Neither $c_1^{BC}({\cal K}_X)$ nor   $-c_1^{BC}({\cal K}_X)$ is
pseudo-effective.

It is  easy to see that in this case, there do exist Gauduchon metrics $g$ for which
${\cal A}$ is stable. Indeed, by Lemma \ref{sigma}, there exist   Gauduchon metrics
$g_-$, $g_0$ on $X$ such that $\deg_{g_-}({\cal K}_X)<0$, $\deg_{g_0}({\cal
K}_X)=0$. If $X$ did not  contain any cycle at all, then by Proposition
\ref{destlb}, any line subbundle of ${\cal A}$ is isomorphic to ${\cal K}_X$, so
stability is guaranteed as soon as $\deg_g({\cal K}_X)<0$. Therefore, ${\cal A}$
will be $g_-$ -stable in this case.

When $X$ contains a single cycle $C$, denote $\nu:=\deg_{g_0}({\cal
O}(C))>0$ and let $\eta$ a closed $(1,1)$-form representing $c_1^{DR}({\cal K}_X)$.
For any sufficiently small $|t|$ the form
$\omega_t:=\omega_{g_0}+t\eta$ is the K\"ahler metric of a 
  Gauduchon metric $g_t$  with
$\deg_{g_t}({\cal K}_X)=-tb_2(X)$, $\deg_{g_t}({\cal O}(-C))=-\nu+t\langle
c_1^{DR}({\cal K}_X),[C]\rangle$. For  sufficiently small positive $t$ one has
$$\deg_{g_t}({\cal O}(-C))<\deg_{g_t}({\cal K}_X)< \frac{\deg_{g_t}({\cal
K}_X)}{2}=\mu_{g_t}({\cal A})<0\ ,
$$
proving that neither ${\cal K}_X$ nor ${\cal O}(-C)$ destabilizes the bundle ${\cal
A}$.

\item $-c_1^{BC}({\cal K}_X)$ is pseudo-effective.

In this case for {\it every} Gauduchon metric $g$ on $X$ the subbundles of ${\cal
A}$ which are isomorphic to ${\cal K}_X$ do not $g$-destabilize ${\cal A}$. 

Therefore if ${\cal A}$ is unstable for every  $g\in{\cal G}(X)$, $X$ must contain
a cycle
$C\subset X$ such that
$$\deg_g({\cal O}(-C))\geq \frac{1}{2}\deg_g({\cal K}_X)\ \forall g\in{\cal G}(X)\ .
$$
\end{enumerate}

By Corollary \ref{equiv} this implies that the Bott-Chern class $c_1^{BC}({\cal
K}_X^\vee (2C))$ is pseudo-effective.
\qed

\begin{co}\label{final} If ${\cal A}$ is unstable for any Gauduchon metric on $X$,
then
$X$ contains a GSS. 
\end{co}
\pf  When $X$ contains two cycles, it must be either a hyperbolic or a parabolic
Inoue surface \cite{Na1}, so it contains a GSS. When $X$ does not contain two
cycles,
$c_1^{BC}(\Theta_X)$ must be pseudo-effective. By Corollary
\ref{curve-ex} it follows that a multiple of the de Rham class $- c_1^{DR}({\cal
K}_X)$ is represented by an effective divisor. In other words, $X$ contains a
numerically pluri-anticanonical divisor. But,  by the main result of \cite{D}, such
a surface contains a GSS.
\qed

Corollary \ref{final} can be reformulated as follows:
\begin{re}
 If $X$ was a counter-example to the GSS conjecture, $X$ must admit
Gauduchon metrics with respect to which the bundle ${\cal A}$  is stable.
\end{re}

The surfaces $X$ with the property ``${\cal A}$ is unstable for every Gauduchon
metric" are very special. Indeed, either $X$ is a (hyperbolic or parabolic) Inoue
surface, or

\begin{re}
If $X$ is in case 2. of Theorem \ref{intth} it must contain a single cycle $C$,
and, writing $c_1^{DR}({\cal O}(C))=-e_I$ with $I\subset\{1,\dots,b_2(X)\}$, a
multiple of the cohomology class
$e_I-e_{\bar I}$ must be represented by an effective divisor.
 
\end{re}
This follows again by  Corollary
\ref{curve-ex} using the pseudo-effectiveness of the class $c_1^{BC}({\cal
K}_X^\vee(-2C))$.
\\ \\
{\bf Acknowledgments:} I am  indebted to Georges Dloussky and Nicholas Buchdahl for
important suggestions and useful  discussions on  the subject. Many  
statements in the present article are direct consequences of their deep 
results.

{\small
Author's address: \vspace{2mm}\\
Andrei Teleman, LATP, CMI,   Universit\'e de Provence,  39  Rue F.
Joliot-Curie, 13453 Marseille Cedex 13, France,  e-mail:
teleman@cmi.univ-mrs.fr. }

\end{document}